\documentclass{article}
\usepackage{amsfonts,amsmath,amssymb}
\usepackage{graphics,epsfig,graphicx}
\usepackage{color}
\allowdisplaybreaks
\begin{document}
\newtheorem{theorem}{Theorem}[section]
\newtheorem{proposition}{Proposition}[section]
\newtheorem{lemma}{Lemma}[section]
\newtheorem{corollary}{Corollary}[section]
\newtheorem{remark}{Remark}[section]
\newtheorem{definition}{Definition}[section]
\renewcommand{\thesection}{\arabic{section}}
\renewcommand{\theequation}{\thesection.\arabic{equation}}
\renewcommand{\thetheorem}{\thesection.\arabic{theorem}}
\numberwithin{equation}{section}
\numberwithin{theorem}{section}
\numberwithin{proposition}{section}
\numberwithin{lemma}{section}
\numberwithin{remark}{section}
\setcounter{secnumdepth}{3}
\newcommand{\cl}{\centerline}
\newcommand{\sms}{\smallskip}
\newcommand{\ms}{\medskip}
\newcommand{\bs}{\bigskip}
\newcommand{\noi}{\noindent}
\newcommand{\itl}[1]{\textit{#1}}
\newcommand{\blf}[1]{\textbf{#1}}
\newcommand{\dsty}{\displaystyle}
\newcommand{\txty}{\textstyle}
\newcommand{\ssty}{\scriptstyle}
\newcommand{\tty}{\texttt}


\newcommand\Par{\mathhexbox278\,}


\newcommand{\al}{\alpha}
\newcommand{\Al}{\Alpha}
\newcommand{\be}{\beta}
\newcommand{\Be}{\Beta}
\newcommand{\Gm}{\Gamma}
\newcommand{\gm}{\gamma}
\newcommand{\dl}{\delta}
\newcommand{\Dl}{\Delta}
\newcommand{\lm}{\lambda}
\newcommand{\Lm}{\Lambda}
\newcommand{\kp}{\kappa}
\newcommand{\varep}{\varepsilon}
\newcommand{\eps}{\epsilon}
\newcommand{\vp}{\varphi}
\newcommand{\sig}{\sigma}
\newcommand{\Sig}{\Sigma}
\newcommand{\om}{\omega}
\newcommand{\Om}{\Omega}
\newcommand{\uom}{\mbox{\boldmath$\omega$}}
\newcommand{\btau}{\mbox{\boldmath$\tau$}}
\newcommand{\bnu}{\mbox{\boldmath$\nu$}}
\newcommand{\up}{\upsilon}
\newcommand{\z}{\zeta}


\newcommand{\df}[1]{\buildrel\mbox{\small def}\over{#1}}
\newcommand{\op}[1]{\buildrel\mbox{\tiny o}\over{#1}}
\newcommand{\db}{\prime\prime}
\newcommand{\bsl}{\backslash}
\newcommand{\lnrm}{\|\!|}
\newcommand{\rnrm}{|\!\|}
\newcommand{\lb}{\lbrack\!\lbrack}
\newcommand{\rb}{\rbrack\!\rbrack}
\newcommand\la{\langle}
\newcommand\ra{\rangle}
\newcommand{\ev}{\equiv}
\newcommand{\nev}{\not\equiv}
\newcommand{\nn}{\mathbb{N}}
\newcommand{\qq}{\mathbb{Q}}
\newcommand{\zz}{\mathbb{Z}}
\newcommand{\rr}{\mathbb{R}}
\newcommand{\rn}{\rr^N}
\newcommand{\cc}{\mathbb{C}}
\newcommand{\id}{\mathbb{I}}
\newcommand{\bo}{\mathbb{O}}

\newcommand{\amsb}[1]{\mathbb{#1}}
\newcommand{\mcl}[1]{\mathcal{#1}}
\newcommand{\bl}[1]{\mathbf{#1}}
\newcommand{\ov}[1]{\overline{#1}}
\newcommand{\wt}[1]{\widetilde{#1}}
\newcommand{\wh}[1]{\widehat{#1}}

\newcommand{\llra}{\leftrightarrow}
\newcommand{\lra}{\longrightarrow}
\newcommand{\LLR}{\Longleftrightarrow}
\newcommand{\LRA}{\Longrightarrow}
\newcommand{\LLA}{\Longleftarrow}


\newcommand{\bbox}{\vrule height.6em width.6em 
depth0em} 
\newcommand{\os}{\vbox{\hrule \hbox{\vrule 
height.6em depth0pt 
\hskip.6em \vrule height.6em depth0em}
\hrule}} 


\newcommand{\Ker}{\operatorname{Ker}}
\newcommand{\Imm}{\operatorname{Im}}
\newcommand{\rank}{\operatorname{rank}}
\newcommand{\dvg}{\operatorname{div}}
\newcommand{\curl}{\operatorname{curl}}
\newcommand{\supp}{\operatorname{supp}}
\newcommand{\essup}{\operatornamewithlimits{ess\,sup}}
\newcommand{\essinf}{\operatornamewithlimits{ess\,inf}}
\newcommand{\essosc}{\operatornamewithlimits{ess\,osc}}
\newcommand{\osc}{\operatornamewithlimits{osc}}
\newcommand{\sign}{\operatorname{sign}}
\newcommand{\loc}{\operatorname{loc}}
\newcommand{\diam}{\operatorname{diam}}
\newcommand{\dist}{\operatorname{dist}}
\newcommand{\card}{\operatorname{card}}
\newcommand{\meas}{\operatorname{meas}}
\newcommand{\spn}{\operatorname{span}}
\newcommand{\dtm}{\operatorname{det}}
%


\newcommand{\overlim}{\mathop{\overline{\lim}}\limits}
\newcommand{\underlim}{\mathop{\underline{\lim}}\limits}
\newcommand{\ttop}[2]{\genfrac{}{}{0pt}{}{#1}{#2}}
\newcommand{\bcu}{\mathop{\txty{\bigcup}}\limits}
\newcommand{\bca}{\mathop{\txty{\bigcap}}\limits}
\newcommand{\bsu}{\mathop{\txty{\sum}}\limits}
\newcommand{\pro}{\mathop{\txty{\prod}}\limits}


\newcommand{\pl}{\partial}
\newcommand{\ptt}{\frac{\pl}{\pl t}}
\newcommand{\ppx}{\frac\pl{\pl x}}
\newcommand{\dds}{\frac d{ds}}
\newcommand{\ddt}{\frac d{dt}}

\newcommand{\intl}{\int\limits}
\newcommand{\iintl}{\iint\limits}
\def\Xint#1{\mathchoice
    {\XXint\displaystyle\textstyle{#1}}%
    {\XXint\textstyle\scriptstyle{#1}}%
    {\XXint\scriptstyle\scriptscriptstyle{#1}}%
    {\XXint\scriptscriptstyle\scriptscriptstyle{#1}}%
    \!\int}
\def\XXint#1#2#3{\setbox0=\hbox{$#1{#2#3}{\int}$}
    \vcenter{\hbox{$#2#3$}}\kern-0.5\wd0}
\def\bint{\Xint-}
\def\dashint{\Xint{\raise4pt\hbox to7pt{\hrulefill}}}
\def\dashiint{\bint\kern-0.15cm\bint}

\newcommand{\ovl}[3]{\int_{#1}^{#2}\kern-#3pt\raise4pt\hbox to7pt{\hrulefill}\ }

\newcommand{\ovll}[3]{\intl_{#1}^{#2}\kern-#3pt\raise4pt\hbox to7pt{\hrulefill}\ }

\newcommand{\tvl}[2]{\iint_{#1}\kern-#2pt\raise4pt\hbox to7pt{\hrulefill}\ }



\newcommand{\omt}{\Om_T}
\newcommand{\plo}{\partial\Omega}
\newcommand{\ovo}{\bar{\Om} }

%
\newcommand{\ci}[1]{C^\infty\!\left({#1}\right)}
\newcommand{\cio}[1]{C_o^\infty\!\left({#1}\right)}
\newcommand{\lloc}[1]{L_{\loc}\!\left({#1}\right)}
\newcommand{\xy}{|x-y|}


\newcommand{\intom}{\intl_{\Om}}
\newcommand{\intbo}{\intl_{\plo}}
\newcommand{\inom}{\int_{\Om}}
\newcommand{\inbo}{\int_{\plo}}
\newcommand{\intrn}{\intl_{\rn}}


\newcommand{\bye}{\end{document}}

\newcommand{\ine}{\int_{E}}
\newcommand{\iinet}{\iint_{E_T}}
\newcommand{\lpf}{\|f\|_{p,E}}
\newcommand{\lqf}{\|f\|_{q,E}}
\newcommand{\ovtau}{\bar{\tau}}
\newcommand{\ovep}{\bar{\eps}}
\newcommand{\ovdl}{\bar{\dl}}
\newcommand{\ovc}{\bar{c}}
\newcommand{\ovg}{\bar{\gm}}
\newcommand{\bry}{B_{\rho}(y)}
\newcommand{\kry}{K_{\rho}(y)}
\newcommand{\ple}{\pl E}
\newcommand{\ove}{\bar{E}}
\newcommand{\trpo}{\Big(\frac{c}{u\pto}\Big)^{p-2}}
\newcommand{\trmo}{\Big(\frac{c}{u\pto}\Big)^{m-1}}
\newcommand{\tpso}{\Big(\frac{u\pto}{c}\Big)^{2-p}}
\newcommand{\datap}{\{p,N,C_o,C_1\}}
\newcommand{\datam}{\{m,N,C_o,C_1\}}
\newcommand{\pto}{(x_o,t_o)}
\def\po{\left(x_o,t_o\right)}
\def\bx{\bar x}
\def\bt{\bar t}
\def\ox{\bar x}
\def\bxt{(\bx,\bt)}
\newcommand{\ukjm}{(u-k_j)_{-}}
\newcommand{\uhm}{(u-h)_{-}}
\newcommand{\umm}{(u-M)_{-}}
\newcommand{\ump}{(u-(\mu_+-M))_+}
\newcommand{\ukpm}{(u-k)_{\pm}}
\newcommand{\ukp}{(u-k)_+}
\newcommand{\uknp}{(u-k_n)_+}
\newcommand{\uknpu}{(u-k_{n+1})_+}
\newcommand{\uknm}{(u-k_n)_-}
\newcommand{\ukm}{(u-k)_-}
\newcommand{\wkp}{(w-k)_+}
\newcommand{\wknp}{(w-k_n)_+}
\newcommand{\wknpu}{(w-k_{n+1})_+}
\newcommand{\wkm}{(w-k)_-}
\newcommand{\vkjm}{(v-k_j)_-}
\newcommand{\rscc}{\frac{e^{\frac{\tau}{p-2}}}{M}
(\dl\rho^p)^{\frac1{p-2}} } 
\newcommand{\ptb}{(\bar{x},\bar{t})}
\newcommand{\vkp}{(v-k)_+}
\newcommand{\vklm}{(v-\bar{\lm}k)_-}
\newcommand{\psiko}{\Psi(H_k^+,\ukp,c)}

\newcommand{\qrtpm}{Q_{\rho}^{\pm}(\theta)}
\newcommand{\qrtp}{Q_{\rho}^+(\theta)}
\newcommand{\qrtm}{Q_{\rho}^-(\theta)}
\newcommand{\qrt}{Q_{\rho}(\theta)}
\newcommand{\qrttm}{Q_{2\rho}^-(\theta)}
\newcommand{\qrttp}{Q_{2\rho}^+(\theta)}
\newcommand{\qrtt}{Q_{2\rho}(\theta)}
\newcommand{\tkn}{\tilde{K}_n}
\newcommand{\trn}{\tilde{\rho}_n}
\newcommand{\trsin}{\left(\frac{e^{\tau}}
{M^{2-p}\dl_1\rho^p}\right)^{\frac1{2-p}}}
\newcommand{\trspo}{\left(\frac{e^{\tau}}
{M^{2-p}\dl_1\rho^p}\right)^{\frac{p-1}{2-p}}}
\newcommand{\ppt}{\frac{\pl}{\pl\tau}}
\newcommand{\kwp}{(k-w)_+}
\newcommand{\kwpo}{[k-(k-w)_++\epsilon k]^{p-1}}
\newcommand{\ikfo}{\int_{K_{8}} }
\newcommand{\enw}{(\eps^n-w)_+ }
\newcommand{\efnp}{\frac{\eps^{n(2-p)}}{[1+\eps-s]^{p-1}}}
\newcommand{\iepnw}{\ikfo\z^p\tau_*\chi_{[\enw>s\eps^n]}dz}
\newcommand{\lngep}{\frac{\gm}{\gm_o}\bigg(
\ln{\frac{1+\eps}{1+\eps-s}}\bigg)^{-p}}
\newcommand{\ipsin}[1]{\ikfo\z^p{#1}
\Psi_{\eps^n}\big[w(z,{#1})\big]dz}
\newcommand{\iphin}[1]{\ikfo\z^p{#1}
\Phi_{\eps^n}\big[w(z,{#1})\big]dz}
\newcommand{\ukno}{(u-k_{n+1})_+}
\newcommand{\tvls}[2]{\iint_{#1}\kern-#2pt\raise4pt\hbox to15pt{\hrulefill}\ }
\newcommand{\uroz}{\int_{K_\rho} u^r(\cdot,0)dx} 
\newcommand{\buroz}{\bint_{K_\rho} u^r(\cdot,0)dx} 
\newcommand{\uqonrm}{\Big(\uqo\Big)^{\frac1q}} 
\newcommand{\urto}{\bint_{K_\rho(x_o)}u^r(\cdot,t_o)dx} 
\newcommand{\uro}{\bint_{K_\rho(x_o)}u^r(\cdot,t_o)dx} 
\newcommand{\uroonrm}{\Big(\urto\Big)^{\frac1r}} 
\newcommand{\uroo}{\bint_{K_{4\rho}(x_o)} 
u^r(\cdot,t_o-\theta_o\rho^2)dx} 
\newcommand{\uronrm}{\Big(\uro\Big^{\frac1r}} 
\newcommand{\Phihn}[1]{\Phi_{h^n}[u(x,{#1})]}
\newcommand{\Psihn}[1]{\Psi^2_{h^n}[u(x,{#1})]}

\def\phi{\varphi}
\def\extr{\Big\arrowvert}
\def\iint{\int\mskip-30mu\int}
\def\Iint{\int\mskip-30mu\int\limits}
\def\mum{\mu_-}
\def\mup{\mu_+}
\def\muduer{\mu(2R)}
\def\mur{\mu(R)}
\def\muq{\mu_{Q(R)}}
\def\umenkm{(u-k)_-}
\def\umenkp{(u-k)_+}
\def\umenpm{(u-k)_{\pm}}
\def\vmenkm{(v-k)_-}
\def\vmenkp{(v-k)_+}
\def\vmenpm{(v-k)_{\pm}}
\def\capparo{K_{\rho}}
\def\capsigro{K_{\sigma\rho}}
\def\niuno{\nu_1}
\def\nidue{\nu_2}
\def\nip{\nu^+}
\def\nim{\nu^-}
\def\csi{\xi}
\def\csip{\xi^+}
\def\csim{\xi^-}
\def\cspm{\xi^{\pm}}
\def\del{\delta}
\def\akrot{A_{k,\rho}}
\def\alrot{A_{l,\rho}}
\def\w1p{W^{1,p}(\rn)}
\def\degpin{[DG]^+_p(\Om)}
\def\degmin{[DG]^-_p(\Om)}
\def\degpmin{[DG]^{\pm}_p(\Om)}
\def\degin{[DG]_p(\Om)}
\def\degpbo{[DG]^+_p(h;\bar\Om)}
\def\degmbo{[DG]^-_p(h;\bar\Om)}
\def\degpmbo{[DG]^{\pm}_p(h;\bar\Om)}
\def\degbo{[DG]_p(h;\bar\Om)}
\def\poten{{\cal U}_K}
\def\media{\mkern12mu\hbox{\vrule height4pt
           depth-3.2pt width5pt}\mkern-16.5mu\int
           \nolimits}
\def\lap{L^{\alpha,p}}
\def\wap{W^{\alpha,p}}
\def\capac{\hbox{cap}}
\newcommand{\dg}{[DG]_p^\pm(E)}
\newcommand{\gdg}{[GDG]_p^\pm(E)}
\newcommand{\dgp}{[DG]_p^+(E)}
\newcommand{\gdgp}{[GDG]_p^+(E)}
\newcommand{\dgm}{[DG]_p^-(E)}
\newcommand{\gdgm}{[GDG]_p^-(E)}
\newcommand{\hk}{\chi_{[k>h]}}

\newenvironment{ack}{\medskip{\it Acknowledgement.}}{}
\title{A Wiener-Type Condition for Boundary Continuity 
of Quasi-Minima of Variational Integrals}
\author{Emmanuele DiBenedetto\footnote{Supported by NSF grant DMS-1265548}\\
Department of Mathematics, Vanderbilt University\\  
1326 Stevenson Center, Nashville TN 37240, USA\\
email: {\tt em.diben@vanderbilt.edu}
\and
Ugo Gianazza\\
Dipartimento di Matematica ``F. Casorati", Universit\`a di Pavia\\ 
via Ferrata 1, 27100 Pavia, Italy\\
email: {\tt gianazza@imati.cnr.it}
}
\date{}
\maketitle
\begin{abstract}
A Wiener-type condition for the continuity at the boundary 
points of Q-minima, is established, in terms of the divergence 
of a suitable Wiener integral [(\ref{Eq:1:7}) and 
Theorem~\ref{Thm:1:1}].
\end{abstract}
\bs
\section{Introduction}\label{S:intro} 
Let $E$ be a bounded, open subset of  $\rn$ and let 
$f:E\times\rr^{N+1}\to\rr$ be a Carath\'eodory 
function satisfying
\begin{equation}\label{Eq:1:1}
C_o |Du|^p\le f(x,u,Du)\le C_1 |Du|^p,
\end{equation}
for constants $0<C_o\le C_1$, and some fixed $p>1$.  
A function $u\in W^{1,p}_{\loc}(\rn)$ 
is a Q-sub(super)minimum for the functional
\begin{equation}\label{Eq:1:2}
J(u)=\ine f(x,u,Du)dx
\end{equation}
if there exists $Q\ge1$ such that 
\begin{equation}\label{Eq:1:3}
J(u)\le Q J\big(u-(+)\vp\big),
\end{equation}
for all non-negative functions $\vp\in W^{1,p}_{\loc}(\rn)$ 
with $\supp\vp\subset\bar{E}$. A function $u\in W^{1,p}_{\loc}(\rn)$ 
is a Q-minimum for $J$ if it satisfies (\ref{Eq:1:3}) 
for all $\vp\in W^{1,p}_{\loc}(\rn)$ with $\supp\vp\subset\bar{E}$ 
and no further sign restriction (\cite{GG}). 
 A Q-minimum $u$ takes boundary values 
$g\in W^{1,p}_{\loc}(\rn)$ on $\ple$ if 
\begin{equation}\label{Eq:1:4}
u-g\in W^{1,p}_o(E).
\end{equation}
If $g\in C(\rn)$ one asks under what conditions on $\ple$, 
the boundary datum $g$ is taken by $u$ 
in the sense of continuous functions. 
 Let $y\in\ple$, denote with $B_\rho(y)$ the ball 
of radius $\rho$  about $y$. For $1<p<N$,
the $p$-capacity of the 
compact set $E^c\cap \bar{B}_\rho(y)$ is defined by
\begin{equation}\label{Eq:1:5}
c_p[E^c\cap\bar{B}_\rho(y)]=\inf_{
\ttop{\psi\in W^{1,p}_o(\rn)\cap C(\rn)}{E^c\cap\bar{B}_\rho(y)
\subset[\psi\ge1]}}\int_{\rn}|D\psi|^pdx.
\end{equation}
For {$1<p< N$}, the relative $p$-capacity of 
$E^c\cap\bar{B}_\rho(y)$ with respect to $B_\rho(y)$ is
\begin{equation}\label{Eq:1:6}
\qquad 
\dl_y(\rho)=\frac{c_p[E^c\cap\bar{B}_\rho(y)]}{\rho^{N-p}}, 
\qquad\qquad (1<p<N).
\end{equation}
If $p=N$, and for $0<\rho<1$, the $N$-capacity of the 
compact set $E^c\cap\bar{B}_\rho(y)$, with respect to the 
ball $B_2(y)$, is defined by 
\begin{equation}\label{Eq:1:5:N}
c_N[E^c\cap\bar{B}_\rho(y)]={\inf_{
\ttop{\psi\in W^{1,N}_o(B_2(y))\cap C_o(B_2(y))}{E^c\cap\bar{B}_\rho(y)
\subset[\psi\ge1]}}\int_{B_2(y)}|D\psi|^Ndx}.
\end{equation}
{The relative capacity $\dl_y(\rho)$ can be formally defined by (\ref{Eq:1:6}), 
for all $1<p\le N$: for $p=N$, $\dl_y(\rho)\equiv c_N[E^c\cap\bar{B}_\rho(y)]$, as defined by (\ref{Eq:1:5:N})}.

For a positive parameter $\eps$ denote by $I_{p,\eps}(y,\rho)$ 
the Wiener integral of $\ple$ at $y\in\ple$, i.e.,   
\begin{equation}\label{Eq:1:7}
I_{p,\eps}(y,\rho)=\int_\rho^1[\dl_y(t)]^{\frac{1}{\eps}}
\frac{dt}{t}.
\end{equation}
The main result of this note is:
\begin{theorem}\label{Thm:1:1}
Let $u$ be a $Q$-minimum for the functional $J(u)$, 
for $1<p\le N$. Assume that $u$  takes a continuous 
datum $u=g$ on $\ple$ in the sense of (\ref{Eq:1:4}). 
There exists $\eps\in(0,1)$, and $\gm>1$, that can be 
determined apriori, quantitatively only in terms 
of $N$, $p$, {$Q$, and the ellipticity ratio $\frac{C_1}{C_o}$}, such that for all $y\in\ple$, 
and all $\rho\in(0,1)$
\begin{equation}\label{Eq:1:8}
\essosc_{E\cap B_\rho(y)} u\le 
\gm\max\Big\{\osc_{\ple\cap B_\rho(y)} g\,;\, 
{\left(\osc_{E\cap B_1(y)} u\right)}\exp\big(-I_{p,\eps}(y,\rho)\big)\Big\}.
\end{equation}
\end{theorem}
Thus, when $1<p\le N$, a Q-minimum $u$, when given 
continuous boundary data $g$ on $\ple$,  is continuos 
up to $y\in\ple$, if the Wiener integral 
$I_{p,\eps}(y,\rho)$ diverges as $\rho\to0$. 
If $p>N$ the continuity of $u$, is insured by the Sobolev 
embedding theorem.
\subsection{Novelty and Significance}\label{S:comments}
The celebrated Wiener criterion states that a harmonic 
function in $E$ is continuous up to $y\in\ple$ if and only 
if the Wiener integral $I_{2,1}(y,\rho)$ diverges as 
$\rho\to0$ (\cite{wiener}). 
Next, for a given $g\in W^{1,p}(\rn)\cap C(\rn)$ consider 
the boundary value problem
\begin{equation}\label{Eq:1:9}
\begin{array}{ll}
u-g\in W^{1,p}_o(E),\quad&\text{for }\>p>1,\\
\dvg{\bf a}(x,u,Du)=0,\quad&\text{weakly in }\> E,
\end{array}
\end{equation}
where, the vector field $\bl{a}$ is subject to the 
structure conditions
\begin{equation}\label{Eq:1:10}
\begin{aligned}
{\bf a}(x,u,Du)\cdot Du&\ge C_o |Du|^p,\\
|{\bf a}(x,u,Du)|&\le C_1|Du|^{p-1}
\end{aligned}
\end{equation}
for constants $0<C_o\le C_1$, and some fixed $p>1$. 
The prototype is
\begin{equation}\label{Eq:1:11}
\begin{array}{ll}
u-g\in W^{1,p}_o(E)\quad&\text{for }\>p>1,\\
\dvg(|Du|^{p-2}Du)=0\quad&\text{weakly in }\> E.
\end{array}
\end{equation}
For solutions of (\ref{Eq:1:11}) Theorem~\ref{Thm:1:1} 
is due to Maz'ja (\cite{maz1}), with the optimal value 
of the parameter $\eps=(p-1)$. 
The proof is based on the comparison principle 
and the Harnack inequality. For solutions of 
(\ref{Eq:1:9})--(\ref{Eq:1:10}) the result is due to 
Gariepy and Ziemer (\cite{gar-zie}), still for optimal 
value of the parameter $\eps=(p-1)$. For these quasi-linear 
equations there is not, in general, a maximum principle. 
Their proof is based on the Moser's logarithmic estimates 
(\cite{moser}) leading to the Harnack inequality for some 
proper convex functions of the solutions, near the boundary 
point $y\in\ple$. In their approach, the structure of 
the p.d.e. in (\ref{Eq:1:9})--(\ref{Eq:1:10}) is crucial. 

Each such quasi-linear equation is the Euler equation of a 
functional $J$, for a suitable integrand $f(x,u,Du)$ (\cite{GG}). 
The notion of Q-minimum is considerably 
more general as it includes {\it almost minimisers}, or 
even minimisers of functionals $J(u)$ which do not 
admit a Euler equation due to the possible lack of 
Gateaux differentiability of $J$.

Nevertheless Q-minima share several crucial properties 
of solutions of quasi-linear equations of the type 
(\ref{Eq:1:9})--(\ref{Eq:1:10}). For example, they are locally 
bounded and locally H\"older continuous in $E$. 
Their interior continuity carries at those boundary 
points where $\ple$ has positive geometric density (\cite{GG}). 
Moreover, non-negative Q-minima satisfy the Harnack inequality 
(\cite{dibe-trud}). However, Q-minima are not known to satisfy 
a maximum principle, nor Harnack inequalities near $\ple$. 

The significance of a Wiener condition for Q-minima, 
is that the structure of $\ple$ near a boundary point 
$y\in\ple$, for $u$ to be continuous up to $y$, hinges 
on minimizing a functional, rather than solving 
an elliptic p.d.e. 


The only result, to date, in this direction, states that 
a Q-minimum $u$, with continuous boundary data $g\in C(\ple)$, 
is continuous up to a boundary point $y\in\ple$ if (\cite{ziem})
\begin{equation}\label{Eq:1:12}
\int_{\rho}^1 \exp\Big(-\frac1{\dl_y(t)^{\frac1{p-1}}}\Big)
\frac{dt}{t}\>\to\>\infty\quad\text{ as }\> \rho\to0.
\end{equation}
Ziemer's proof follows from a standard DeGiorgi iteration 
technique (\cite{DeG}). The novelty of our Theorem~\ref{Thm:1:1} is 
in replacing the exponential decay (\ref{Eq:1:12}) 
in the Wiener integral with a power-like decay. 
The technical novelty is in extending a weak Harnack 
inequality for quasi minima (\cite{dibe-trud}), to hold 
near the boundary, coupled with proper choices of 
test functions in (\ref{Eq:1:3}) as indicated by 
Tolksdorf (\cite{TK}).  The optimal value of the 
parameter $\eps=(p-1)$, remains  elusive. 
\vskip.2truecm
\begin{ack}
We  are grateful to the anonymous referee for the valuable comments. After submitting this note, we learnt that a similar result has been proved with a different technique by J. Bj\"orn (see \cite{bjorn}).
\end{ack}
\section{Main Tools in the Proof of Theorem~\ref{Thm:1:1}}\label{S:main_tools}
\subsection{$Q$-Subminima and Test Functions}\label{S:test}
\begin{proposition}\label{Prop:2:1}
Let $y\in \ple$ and let $u$ be a non-negative $Q-subminimum$ 
for $J$, in $\bar{B}_{\rho}(y)\cap\ove$, such that $u=0$ on 
$B_\rho(y)\cap\ple$. There is a positive constant 
$\gm_o$ that can be determined apriori only in terms of 
$N$, $p$, $Q$, {and the ellipticity ratio $\frac{C_1}{C_o}$}, such that
\begin{equation}\label{Eq:2:1}
\int_{B_\rho(y)\cap E}|Du|^p|\vp|^pdx\le
\gm_o\int_{B_\rho(y)\cap E} u^p|D\vp|^pdx,
\end{equation}
for all non-negative $\vp\in W_o^{1,p}\big(B_\rho(y)\big)$.
\end{proposition}
Note that $\vp$ is not required to vanish on $B_\rho(y)\cap\ple$. 
The proof results from a minor variant of an argument of 
Tolksdorf \cite{TK}. From the property (\ref{Eq:1:1}) of 
$f$ and the definition (\ref{Eq:1:2})--(\ref{Eq:1:3}) 
of Q-subminimum, 
\begin{equation}\label{Eq:2:2}
\int_{B_\rho(y)\cap E}|Du|^p dx\le
Q\frac {C_1}{C_o}\int_{B_\rho(y)\cap E} |D(u-u\vp)|^pdx,
\end{equation}
for all non-negative $\vp\in W_o^{1,p}\big(B_\rho(y)\big)$. 
The new observation here is that since $u$ vanishes on 
$B_\rho(y)\cap\ple$, the test function $u\vp$ is admissible 
in (\ref{Eq:1:3}) even if $\vp$ does not vanish on 
$B_\rho(y)\cap\ple$, provided it does vanish on 
$\pl B_\rho(y)$. The remaining arguments leading to 
(\ref{Eq:2:1}) starting from (\ref{Eq:2:2}) are identical 
to those in \cite{TK}.
\begin{corollary}\label{Cor:2:1}
Let $u$ satisfy the same assumptions as 
Proposition~\ref{Prop:2:1}. Then for all constants $h>0$ 
\begin{equation}\label{Eq:2:3}
\int_{B_\rho(y)\cap E}|D(u+h)|^p|\vp|^pdx\le
\gm_o\int_{B_\rho(y)\cap E} (u+h)^p|D\vp|^pdx,
\end{equation}
for all non-negative $\vp\in W_o^{1,p}\big(B_\rho(y)\big)$.
The constant $\gm_o$ is the same as in (\ref{Eq:2:1}) 
and is independent of $h$. 
\end{corollary}
\subsection{$Q$-Superminima and the Weak Harnack 
Inequality}\label{S:weak_harnack}
\begin{proposition}\label{Prop:2:2}
Let $y\in \ple$ and let $v\in W^{1,p}(B_{2\rho}(y))$ 
be non-negative and satisfying 
\begin{equation}\label{Eq:2:4}
\int_{B_r(z)}|D(v-k)_-|^p dx\le{\frac{\gm_1}{r^p}}\int_{B_{2r(z)}} (v-k)_-^pdx
\end{equation}
for all balls $B_{2r}(z)\subset B_{2\rho}(y)$ and all $k>0$, for 
a constant $\gm_1$ independent of $k$, $z$ and $r$. Then, 
there exist constants $C>1$ and $\eps\in(0,1)$, that can 
be determined apriori only in terms of $N$, $p$, and the 
{constant $\gm_1$ in (\ref{Eq:2:4})}, 
such that 
\begin{equation}\label{Eq:2:5}
\left(\frac1{|B_\rho(y)|}\int_{B_\rho(y)}v^{\eps}dx
\right)^{\frac1{\eps}}\le C 
\essinf_{B_\rho(y)} v.
\end{equation}
\end{proposition}
The weak Harnack inequality (\ref{Eq:2:5}) is a sole 
consequence of the family of inequalities (\ref{Eq:2:4}), 
and as such, disconnected from the notion of Q-superminimum 
(\cite{dibe-trud}).  However, if $v$ is a Q-superminimum 
in $E$, for balls $B_{2\rho}(y)\subset E$, inequalities 
(\ref{Eq:2:4}) are satisfied by $v$ (\cite{GG}). 
\section{Proof of Theorem~\ref{Thm:1:1}}\label{S:proof_thm}
\subsection{Estimating the Oscillation About a Point 
$y\in\ple$ by the 
Weak Harnack Inequality}\label{S:proof_thm_001}
Having fixed $y\in\ple$ assume without loss of generality 
that $y=0$ and write $B_\rho(0)=B_\rho$, and continue to 
denote by $g$ the boundary datum of $u$, in the sense of 
(\ref{Eq:1:4}). We may assume that at least one of the 
following two inequalities holds true:
\begin{align*}
\essup_{B_{2\rho}\cap E}u-{\txty\frac14} \essosc_{B_{2\rho}\cap E}u
&>\essup_{B_{2\rho}\cap\ple}g;\\
\essinf_{B_{2\rho}\cap E}u\,{+}\,{\txty\frac14} \essosc_{B_{2\rho}\cap E}u
&<\essinf_{B_{2\rho}\cap\ple}g.
\end{align*}
Indeed, if both are violated one has
\begin{equation*}
\essosc_{B_{2\rho}\cap E}u\le 2 \essosc_{B_{2\rho}\cap\ple}g,
\end{equation*}
and the assertion of the theorem follows. Assuming then 
that the first holds, the function 
\begin{equation*}
\Big(u\,-\big(\essup_{B_{2\rho}\cap E}u-{\txty\frac14}
\essosc_{B_{2\rho}\cap E} u\big)-(1-k){\txty\frac14}
\essosc_{B_{2\rho}\cap E} u\Big)_+ 
\end{equation*}
is a non-negative Q-subminimum, for $J$, in 
$\bar{B}_{2\rho}\cap\ove$, for all $0<k\le1$, 
vanishing on $B_{2\rho}\cap\ple$. As such it satisfies 
(\ref{Eq:2:1}) of Proposition~\ref{Prop:2:1}, 
over $B_{2\rho}$, which we rewrite as 
\begin{equation}\label{Eq:3:1}
\begin{aligned}
\int_{B_{2\rho}\cap E}&|D\big(w-(1-k)\big)_+|^p|\vp|^pdx\\
&\le\gm_o\int_{B_{2\rho}\cap E} \big(w-(1-k)\big)_+^p|D\vp|^pdx,
\end{aligned}
\end{equation}
for all non-negative $\vp\in W_o^{1,p}(B_{2\rho})$, where 
\begin{equation*}
w\df{=}\frac{{\dsty \Big(u\,-\big(\essup_{B_{2\rho}\cap E}u
-{\txty\frac14}\essosc_{B_{2\rho}\cap E}u\big)\Big)_+}}{{\txty\frac14}
{\dsty\essosc_{B_{2\rho}\cap E}u}},
\end{equation*}
for all $0<k\le1$. From the definitions one verifies 
that $0\le w\le1$, and it vanishes on $B_{2\rho}\cap \ple$. 
We continue to denote by $w$ and $\big(w-(1-k)\big)_+$ 
their extensions with zero on $B_{2\rho}\cap E^c$. 
By Corollary~\ref{Eq:2:1}, inequalities (\ref{Eq:3:1}) 
continue to hold for all $k\ge0$. Set $v=1-w$ and 
rewrite (\ref{Eq:3:1}) in the form
\begin{equation}\label{Eq:3:2}
\int_{B_{2\rho}}|D(v-k)_-|^p|\vp|^pdx\le
\gm_o\int_{B_{2\rho}}(v-k)_-^p|D\vp|^pdx,
\end{equation}
for all non-negative $\vp\in W_o^{1,p}(B_{2\rho})$, 
and for all $k\ge0$. In what follows we denote by $\gm$ 
a generic, positive constant that can be quantitatively 
determined apriori only in terms of $N$, $p$, $Q$. 

For a ball $B_{2r}(z)\subset B_{2\rho}$, in (\ref{Eq:3:2}) 
choose $\vp$ as the standard, non-negative cutoff function 
in $B_{2r}(z)$ which equals 1 on $B_{r}(z)$ and such that 
$|D\vp|\le r^{-1}$. For such a choice $(v-k)_-$ satisfies, the 
assumptions of Proposition~\ref{Prop:2:2}. Hence there 
exists $\gm>1$ and $\eps\in(0,1)$ that can be determined 
apriori only in terms of $N$, $p$, {$Q$, and the ellipticity ratio $\frac{C_1}{C_o}$},
such that
\begin{equation}\label{Eq:3:3}
\dashint_{B_{2\rho}} v^{\eps} dx\le \gm^{\eps} 
\left(\frac{\dsty \essup_{B_{2\rho}\cap E} u-
\essup_{B_{\rho}\cap E} u}{\frac14{\dsty
\essosc_{B_{2\rho}\cap E}u}}\right)^{\eps} 
\end{equation}
\begin{remark}\label{Rmk:3:1} {\normalfont
Whence the parameter $\eps$ has been identified, inequality 
(\ref{Eq:3:3}) continues to hold for smaller $\eps$, 
with the same constant $\gm$.
}
\end{remark} 
\subsection{Estimating the Oscillation About a Point 
$y\in \ple$ by the Capacity of 
$E^c\cap\bar{B}_{\rho}(y)$}\label{S:proof_thm_002}
Continue to assume $y=0$ and write $B_\rho(0)=B_\rho$.
\begin{proposition}\label{Prop:3:1}
There exists $p_o\in(1,p)$, that depends only on $N$, $p$, {$Q$, and the ellipticity ratio $\frac{C_1}{C_o}$}, such that for all $p_o\le q<p$, 
and for all non-negative $\z\in W_o^{1,p}(B_{2\rho})$,
there holds 
\begin{equation}\label{Eq:3:4}
\int_{B_{2\rho}} v^{-q}|Dv|^p\z^pdx\le \gm 
\int_{B_{2\rho}} v^{p-q}|D\z|^pdx
\end{equation}
for a constant $\gm>1$ that depends only on $N$, $p$, $Q$, 
$q$, $p_o$, {and the ellipticity ratio $\frac{C_1}{C_o}$}.
\end{proposition} 
\noi{\bf Proof:} Using an idea of \cite{mi-zi}, set 
$\vp=v^\sig\z$ in (\ref{Eq:3:2}) where $\sig\in(0,1)$ 
is a parameter to be chosen and $\z\in W_o^{1,p}(B_{2\rho})$ 
is non-negative. For such choices (\ref{Eq:3:2}) yields
\begin{equation*}
\int_{B_{2\rho}}|D(v-k)_-|^p v^{\sig p}\z^p dx
\le\gm\int_{B_{2\rho}}(v-k)_-^p
\left[\sig^p v^{(\sig-1)p}|Dv|^p\z^p
+v^{\sig p}|D\z|^p\right]dx.
\end{equation*}
Choose $\sig>0$ and $1<q<p$ so that $(1-\sig)p<q$, multiply 
both sides of this inequality by $k^{-\sig p-q-1}$ and integrate 
in $dk$ over $(0,\infty)$. Interchanging the order of integration 
with the aid of Fubini's theorem, the left-hand side equals 
\begin{equation*}
\int_0^{\infty}\int_{B_{2\rho}}|D(v-k)_-|^p
v^{\sig p}\z^pk^{-\sig p-q-1}dxdk=
\frac{1}{\sig p+q}\int_{B_{2\rho}}|Dv|^p v^{-q}\z^pdx.
\end{equation*}
The right-hand side is transformed and estimated by
\begin{align*}
\int_0^{\infty}&\int_{B_{2\rho}}(v-k)_-^p
\left[\sig^p v^{(\sig-1)p}|Dv|^p\z^p+v^{\sig p}|D\z|^p\right]
k^{-\sig p-q-1}dxdk\\
&=\sig^p\int_{B_{2\rho}}v^{(\sig-1)p}|Dv|^p\z^p \left(
\int_v^{\infty}k^{(1-\sig)p-q-1}dk\right)dx\\
\quad &+\int_{B_{2\rho}} v^{\sig p} |D\z|^p\left(
\int_v^{\infty}k^{(1-\sig)p-q-1}dk\right)dx\\
&=\frac{1}{q-(1-\sig)p}\int_{B_{2\rho}}v^{(1-\sig)p-q}
\left[\sig^p v^{-(1-\sig)p}|Dv|^p\z^p+v^{\sig p}|D\z|^p\right]dx.
\end{align*}
Combining these estimates yields
\begin{align*}
\int_{B_{2\rho}}v^{-q}|Dv|^p\z^pdx\le 
&\gm\frac{\sig p+q}{q-(1-\sig)p}\sig^p\int_{B_{2\rho}}
v^{-q}|Dv|^p\z^p\\
&+\gm\frac{\sig p+q}{q-(1-\sig)p}\int_{B_{2\rho}}
v^{p-q}|D\z|^p dx.
\end{align*}
To conclude the proof choose $\sig\in(0,1)$ such that
\begin{align*}
\gm\frac{\sig p+q}{q-(1-\sig)p}\sig^p=\frac12,\quad\text{ and }\quad 
(1-\sig)p<q<p. 
\end{align*}
One may first choose $p_o=(1-\sig^2)p\le q<p$ and then 
$\sig$ so small that the first of these inequalities 
is in force.\hfill\bbox

We now conclude the proof of the Theorem, still following 
\cite{mi-zi}. Fix $p_o\le q<p$ where $p_o$ is the parameter  
claimed in Proposition~\ref{Prop:3:1}, and rewrite it 
as $q=p-\eps$. By virtue of Remark~\ref{Rmk:3:1}, this value of 
$\eps$ can be taken equal to the analogous in (\ref{Eq:3:3}). 
For such a choice, (\ref{Eq:3:4}) gives
\begin{equation*}
\int_{B_{2\rho}} |D[v^{\frac\eps p}\vp]|^p\,dx\le \gm(\eps)
\int_{B_{2\rho}}v^\eps |D\vp|^pdx.
\end{equation*}
Next choose $\vp\in W_o^{1,p}(B_{2\rho})$ to be the standard, 
non-negative cutoff function in $B_{2\rho}$ which equals 1 
on $B_{\rho}$ and such that $|D\vp|\le \rho^{-1}$. For such a 
choice and $\rho$ sufficiently small $v\vp=1$ on $B_\rho\cap E^c$ 
and therefore,
\begin{equation*}
c_p\big[E^c\cap\bar{B}_\rho] \le \frac{\gm(\eps)}{\rho^p}
\int_{B_{2\rho}}v^\eps dx. 
\end{equation*}
Dividing by $\rho^{N-p}$ and combining the resulting 
inequality with (\ref{Eq:3:3}) gives 
\begin{equation*}
\dl^{\frac1{\eps}}_o(\rho)\le\gm\ 
\frac{\dsty \essup_{B_{2\rho}\cap E} u-
\essup_{B_{\rho}\cap E} u}{\frac14{\dsty
\essosc_{B_{2\rho}\cap E}u}}. 
\end{equation*}
This in turn implies
\begin{equation*}
\essosc_{B_\rho\cap E} u\le \Big(1
-\frac1{4\gm}\dl_o(\rho)\Big)\essosc_{B_{2\rho}} u
\end{equation*}
Iteration of this inequality over a sequence of balls 
of dyadic radii $\rho_{-n}=2^{-n}\rho$ yields the Theorem. 


\bye